\documentclass[preprint,12pt,3p,authoryear,nopreprintline]{elsarticle}

\usepackage{amssymb}
\usepackage{amsmath}
\usepackage{amsthm}

\usepackage{graphicx}      
\usepackage{natbib}        
\usepackage{array}
\usepackage[caption=false,font=normalsize,labelfont=sf,textfont=sf]{subfig}
\usepackage{amsmath,amssymb,amsfonts}
\usepackage{mathtools}
\usepackage{graphicx}
\usepackage{xcolor}
\usepackage{tikz}
\usepackage{pgfplots}
\usepackage{comment}
\usepackage{placeins}
\usepackage{booktabs}
\usepackage{bm}
\usepackage{acro}
\pgfplotsset{width=10cm,compat=1.18}
\usetikzlibrary{intersections}
\usepackage{algorithm}
\usepackage{algpseudocode}
\usepackage{hyperref}
\usepackage{textcomp}
\usepackage{nicefrac}

\newcommand{\publicationnotice}{%
  Published in \textit{Engineering Applications of Artificial Intelligence},
  Vol.~181, Part~6, Article~115698, Elsevier, 2026.
  DOI: \href{https://doi.org/10.1016/j.engappai.2026.115698}%
  {10.1016/j.engappai.2026.115698}.%
}
\newcommand{\licensenotice}{%
  \textcopyright~2026. This manuscript version is made available under the CC-BY-NC-ND 4.0 license
  \url{https://creativecommons.org/licenses/by-nc-nd/4.0/}.%
}
\usepackage{eso-pic}
\AddToShipoutPictureBG*{%
  \AtPageUpperLeft{%
    \raisebox{-25mm}[0pt][0pt]{%
      \makebox[\paperwidth][c]{%
        \parbox[t]{\textwidth}{\centering\scriptsize\publicationnotice}%
      }%
    }%
  }%
  \AtPageLowerLeft{%
    \raisebox{15mm}[0pt][0pt]{%
      \makebox[\paperwidth][c]{%
        \parbox[b]{\textwidth}{\scriptsize\raggedright\licensenotice}%
      }%
    }%
  }%
}

\makeatletter
\def\namedlabel#1#2{\begingroup
    #2%
    \def\@currentlabel{#2}%
    \phantomsection\label{#1}\endgroup
}
\makeatother
\algdef{SNb}[MYBLOCK]{Output}{EndOutput}{\textbf{output }}{}
\algdef{SNb}[MYBLOCK]{Input}{EndInput}{\textbf{Input: }}{}

\newtheorem{assumption}{Assumption}

\newtheorem{theorem}{Theorem}

\newtheorem{proposition}{Proposition}
\newtheorem{remark}{Remark}

\DeclareMathOperator*{\argmin}{arg\,min}

\newcommand{\EV}{\widehat{V}}

\newcommand{\innerprod}[1]{{\langle #1\rangle}}

\newcommand{\concat}[2]{\begin{bmatrix}
  #1\\
  #2
\end{bmatrix}}
\newcommand{\KME}[1]{\mathcal{E}_{#1}}
\newcommand{\EKME}[1]{\widehat{\mathcal{E}}_{#1}}
\newcommand{\defeq}{\coloneqq}
\newcommand{\tr}{\intercal}
\newcommand{\order}[1]{{\mathcal{O}}\left(#1\right)}

\DeclareAcronym{rmse}{
  short = RMSE,
  long  = Root Mean Squared Error,
}

\newcommand{\SFA}{{\mathsf{A}}}
\newcommand{\SFB}{{\mathsf{B}}}

\newcommand{\SFX}{\mathsf{X}}

\newcommand{\SFZ}{\mathsf{Z}}
\newcommand{\SFU}{\mathsf{U}}
\newcommand{\SFV}{\bm{\mathsf{v}}}
\newcommand\Set[1]{\mathbb{#1}} 

\newcommand{\SFD}{\bm{D}}
\newcommand{\Kx}{{\bm{K}}_{\mathsf{X}}}
\newcommand{\bxi}{{\bm{x}^{(i)}}}

\newcommand{\bui}{{\bm{u}^{(i)}}}

\newcommand{\bu}{{\bm{u}}}
\newcommand{\lilsum}{\textstyle\sum}
\newcommand{\IN}{{{\mathbb X}}}
\newcommand{\OUT}{{{\mathbb Z}}}
\newcommand{\spIN}{{\mathcal{H}_{\IN}}}
\newcommand{\spOUT}{{\mathcal{H}_{\OUT}}}
\newcommand{\bz}{{\bm{z}}}
\newcommand{\featx}[1]{k(\cdot, #1)}
\newcommand{\bx}{{\bm{x}}}
\newcommand{\stgcost}{{\ell}}

\usepackage{pifont}

\usepackage{xcolor}

\makeatletter
\newcommand{\ccell}[3][]{%
  \kern-\fboxsep
  \if\relax\detokenize{#1}\relax
    \expandafter\@firstoftwo
  \else
    \expandafter\@secondoftwo
  \fi
  {\colorbox{#2}}%
  {\colorbox[#1]{#2}}%
  {#3}\kern-\fboxsep
}
\makeatother
\definecolor{cellclr}{gray}{0.9}
\AtBeginDocument{\setlength{\cmidrulekern}{0.3em}}

\newcommand{\stmbxi}{{\boldsymbol{x}^{(i)}}}

\begin{document}
\begin{frontmatter}
\title{Data-Driven Stochastic Optimal Control in Reproducing Kernel Hilbert Spaces \tnoteref{thanks}}
\tnotetext[thanks]{This work was supported by the DAAD program Konrad Zuse Schools of Excellence in Artificial Intelligence, sponsored by the Federal Ministry of Education and Research, and by the European Union's Horizon Europe innovation action program under grant agreement No. 101093822, ``SeaClear2.0''.}

\author[First]{Nicolas Hoischen}
\author[First]{Petar Bevanda}
\author[First]{Stefan Sosnowski}
\author[First]{Sandra Hirche}
\author[Second]{Boris Houska}
\nonumnote{\raggedright\textit{Email addresses:}
  \texttt{\char123 nicolas.hoischen, petar.bevanda, sosnowski, hirche\char125 @tum.de}
  \texttt{borish@shanghaitech.edu.cn}.}

\affiliation[First]{organization={Chair of Information-oriented Control, Technical University of Munich},
  country={Germany}}
\affiliation[Second]{organization={School of Information Science and Technology, ShanghaiTech University},
  country={China}}

\begin{abstract}
This paper proposes a fully data-driven approach for optimal control of nonlinear control-affine systems represented by a stochastic diffusion. The focus is on the scenario where both the nonlinear dynamics and stage cost functions are unknown, while only a control penalty function and constraints are provided. To this end, we embed state probability densities into a reproducing kernel Hilbert space (RKHS) to leverage recent advances in operator regression, thereby identifying Markov transition operators associated with controlled diffusion processes. This operator learning approach integrates naturally with convex operator-theoretic Hamilton--Jacobi--Bellman recursions that scale linearly with state dimensionality, effectively solving a wide range of nonlinear optimal control problems. Numerical results demonstrate its ability to address diverse nonlinear control tasks, including the depth regulation of an autonomous underwater vehicle.
\end{abstract}

\end{frontmatter}

\section{Introduction}

Data-driven control methods offer a compelling alternative when deriving accurate first-principle models becomes prohibitively complex. Classical approaches retain structural knowledge of the system dynamics while using
data to estimate unknown parameters through system identification and adaptive control techniques \citep{astrom1994adaptive}. As system complexity grows, methods relying less on prior knowledge have emerged: Gaussian
processes learn unknown dynamics with quantified uncertainty \citep{hewing2019cautious}, data-enabled predictive control synthesizes controllers directly from input--output data \allowbreak\mbox{\citep{coulson2019data}}, and reinforcement
learning methods such as \citet{schulman2017proximal} learn optimal policies through direct environmental interaction, surpassing traditional model-based optimal control---for example, in agile drone racing
\citep{song2023reaching}. This spectrum, from parameter estimation to purely data-driven policies, reflects a fundamental shift in which data becomes the primary source of insight as complexity increases.

In recent years, operator-theoretic approaches have garnered substantial interest for learning nonlinear system dynamics from data \citep{Kostic2022LearningSpaces,Klus2020eig}. Unlike conventional nonlinear representations
in the physical state space, these methods describe the dynamics through a linear operator, most notably, the Koopman operator~\citep{Koopman1931}, which fully characterizes the system's evolution. This representational
simplicity has enabled diverse applications including system identification \citep{li2017extended}, soft robotics \citep{bruder2020data}, robust control \citep{STRASSER20232257} and optimal control
~\citep{Vaidya2024,Caldarelli2024,Korda2018}.

Optimal control in state space traditionally relies on Bellman's principle of optimality, but is severely limited by the curse of dimensionality, which renders Hamilton--Jacobi Bellman equations intractable even for
moderate-scale problems \citep{bertsekas1996stochastic}. In such scenarios, local optimal control approaches, such as direct methods combined with sequential quadratic programming or interior points
solvers~\citep{Diehl2009} have demonstrated superior performance over dynamic programming approaches. However, these local approaches require accurate initialization, are challenged by nonconvexity, and may result in
suboptimal local solutions.

In contrast, convex reformulations offer a promising route to globally optimal solutions in nonlinear control via linear embedding methods \citep{Vinter1993}, sum-of-squares dual approaches \citep{henrion2008nonlinear}, or
operator-theoretic system representations \citep{houska2025convex,guo2022tutorial}, each entailing fundamentally different algorithmic frameworks. Moreover,
\emph{convex} reformulations circumvent the linearization trade-offs inherent to
traditional local methods. Operator-theoretic representations are particularly attractive as they can be learned from a single dataset by modern nonparametric operator-theoretic learning
methods~\citep{Kostic2022LearningSpaces,Kostic2023KoopmanEigenvalues,bevanda2025nonparametric} and integrate naturally with stochastic control theory for diffusion processes \citep{Fleming1993,Bogachev2015,Krylov2008}. We
therefore build on these techniques to address stochastic optimal control design in flexible, infinite-dimensional, Reproducing Kernel Hilbert Spaces (RKHSs), which have a rich history in machine learning
\citep{IngoSteinwart2008SupportMachines}.

\emph{Contributions}: This paper proposes a novel algorithm for solving data-driven stochastic optimal control problems. Our main contributions are:
\begin{enumerate}
\item We develop a nonparametric framework to learn finite-rank approximations of the Markov transition operators of stochastic control systems by using kernel regression techniques. The sample and time complexity of our kernel-based system identification method scales \emph{linearly} with the number of states and controls.
\item We propose an RKHS embedded dynamic programming recursion, which leverages the data-driven Markov operator estimates to compute an associated approximation of the global minimizer of stochastic optimal control problems.

\end{enumerate}
Unlike policy-iterative approaches \citep{meng2024physics} or approaches whose practical implementation relies on blending local pre-stabilizing controllers with global ones, thus losing exact optimality near the
origin \citep{Vaidya2024}, our RKHS-based formulation extracts the optimal policy directly from the value function without an initial stabilizing policy.
Whereas Extended Dynamic Mode Decomposition (EDMD) and related Koopman approaches first choose a dictionary of observables and then identify a finite-dimensional operator on the span of that dictionary
\citep{williams2015data,klus2016numerical,STRASSER20232257,nuske2023finite,bruder2020data}, our nonparametric approach does not require the explicit construction of a finite-dimensional basis dictionary. Instead, the
approximation space is induced by the kernel and the sampled data, so that the effective basis adapts to the data. By utilizing suitable infinite-dimensional RKHSs, it can approximate arbitrary operators to arbitrary
accuracy \citep{Kostic2023KoopmanEigenvalues,IngoSteinwart2008SupportMachines}. Finally, our framework naturally accommodates general non-quadratic stage costs in contrast to recent kernel-based linear quadratic regulator
(LQR) analogs \citep{Caldarelli2024}, and circumvents sum-of-squares synthesis bottlenecks \citep{moyalan2023data}, thereby relaxing common requirements on initial feasible policies, explicit dictionary design, and
restrictive cost structures. In contrast to \citet{bevanda2025kernel}, which relies on specific control excitations to learn infinitesimal generators, our method requires no such conditions, making it applicable to
real-system data.

\emph{Organization}:
Section~\ref{sec:probstat} reviews methods for reformulating stochastic optimal control problems as equivalent Partial Differential Equation (PDE) constrained optimization problems.  Section~\ref{sec:optctrlKME} leverages empirical embeddings in an RKHS to introduce a data-driven approximation of the aforementioned optimal control problem. Finally,  Section~\ref{sec:Numerical} applies the proposed method to benchmark optimal control problems.

\emph{Notation}:
We write $[n] \coloneqq \{1,\dots,n\}$ for integers and $\odot$ for the Hadamard product, $\otimes$ for the tensor product; $\IN = \mathbb{R}^{n_x}$ denotes the state space and $\mathbb U = \mathbb{R}^{n_u}$ the control space. The set of $k$-times continuously differentiable functions is denoted by $C^k(\IN)$, and $L^k(\IN)$ denotes the set of $L^k$-integrable functions. $M_+(\mathbb{X})$ denotes the set of non-negative bounded Borel measures on $\mathbb{X}$. Additionally, $H^k(\IN)$ denotes the Hilbert-Sobolev space of $k$-times weakly differentiable functions on $\mathbb X$, including the special case $L^2(\mathbb X) = H^0(\mathbb X)$. Generally, the dual Hilbert space of a Hilbert space $H$ is denoted by $H^{*}$. As usual, $W(0,T; L^2(\IN), H^1(\IN))$ denotes the Bochner-Sobolev space of functions $\phi \in L^2(0,T; H^1(\IN))$ whose weak time derivatives exist and satisfy $\dot{\phi} \in L^2(0,T; H^1(\IN)^{*})$. With $\mathcal F$ and $\mathcal Y$ denoting separable Hilbert spaces, $\mathrm{HS}({\mathcal{F},\mathcal{Y}})$ denotes the space of Hilbert--Schmidt operators with norm $\| A\| _{\mathrm{HS}({\mathcal{F},\mathcal{Y}})}^2 \defeq \sum_i \| Ae_i\| _{\mathcal{Y}}^2$ for any orthonormal basis $\{e_i\}$ of $\mathcal{F}$. Moreover, ${\mathcal Y} \otimes {\mathcal F}$ denotes the tensor product of Hilbert spaces. Whenever it is clear from the context which Hilbert space is meant, we write $\innerprod{\cdot, \cdot}$ for the inner product. Finally, ``div'' denotes the divergence and ``$\Delta$'' the Laplacian.

\subsection{Informal Problem Statement}\label{sec:informalProblem}

Our goal is to find a \textit{feedback control law} $\boldsymbol{\pi}$ solving the stochastic optimal control problem  on the domain $\mathbb{X} = \mathbb R^{n_x}$
\begin{subequations}
\label{eq::SOCP}
\begin{align}
&\lim_{T \to \infty}  \,  \underset{\boldsymbol{\pi} \in \Pi}{\min}\   \mathbb{E} \left[ \frac{1}{T} \int_0^T \left( \stgcost(\boldsymbol{X}_t) + r( \boldsymbol{\pi}(\boldsymbol{X}_t) ) \right) \, \mathrm{d}t \right] \label{eq::ErgodicCost}\\
& \text{s.t.}
\mathrm{d}\boldsymbol{X}_t \! = \! \left( \boldsymbol{f}(\boldsymbol{X}_t)+\boldsymbol{G}(\boldsymbol{X}_t) \boldsymbol{\pi}( \boldsymbol{X}_t ) \right) \mathrm{d}t + \sqrt{2 \epsilon}~\mathrm{d} \boldsymbol{W}_t, \label{eq::SDE}
\end{align}
\end{subequations}
for \textit{unknown} functions $(\boldsymbol{f},\boldsymbol{G},\ell)$ based on sampled measurements of the system states, inputs, and stage costs, $\{(\boldsymbol{x}_k,\boldsymbol{u}_k), \boldsymbol{x}_{k+1}, \ell(\bx_k)\}$, corresponding to discrete-time observations of the Stochastic Differential Equation (SDE) \eqref{eq::SDE} under a small sampling time $h > 0$.
Moreover, in this context, $T > 0$ denotes the time horizon, $\boldsymbol{W}_t$, $t \geq 0$, denotes a Wiener process, $\epsilon > 0$ a diffusion parameter, and $r$ denotes the control penalty. For brevity, we further introduce the notation
\begin{equation*}
    \Pi \defeq \{ \boldsymbol{\pi} \in L^\infty(\mathbb{X}, \mathbb{R}^{n_u}) \mid \boldsymbol{\pi} \in \mathcal U ~\text{on}~\mathbb{X} \}
\end{equation*}
where the set $\mathcal U  \subseteq \mathbb R^{n_u}$ models the control constraints.

We aim at a nonparametric and fully data-driven method that does not rely on prior knowledge about the system's properties and that
\begin{enumerate}
\item learns \eqref{eq::SDE} and the stage cost from snapshot data;
\item computes a data-driven approximation of the global minimizer of \eqref{eq::SOCP} in a computationally simple manner.
\end{enumerate}

\section{Stochastic Optimal Control}\label{sec:probstat}
In this section, we introduce the stochastic optimal control problem formulation and briefly review select existing findings on the use of Fokker--Planck-Kolmogorov (FPK) operators for solving such problems through
PDE-constrained optimization techniques \citep{Crandall1992,houska2025convex}.
\begin{assumption}\label{ass::blanket}
Throughout this paper, the following blanket assumptions are made.
\begin{enumerate}
\item[\namedlabel{itm:1con}{(1)}] The functions {$\boldsymbol{f} \in C^1(\mathbb{X})^{n_x}$} and $\boldsymbol{G} \in C^1(\mathbb{X})^{n_x \times n_u}$ are continuously differentiable; and $\epsilon >0$ is given.
\item[\namedlabel{itm:2con}{(2)}] The stage cost $\stgcost\in C^1(\mathbb{X})$ is continuously differentiable, bounded from below, and radially unbounded.
\item[\namedlabel{itm:3con}{(3)}] The initial state {$\boldsymbol{X}_0$} is a random variable with given probability distribution $\rho_0 \in L^2(\mathbb{X})$.
\item[\namedlabel{itm:4con}{(4)}] The set {$\mathcal U \subseteq \mathbb R^{n_u}$} is closed, convex, and $\mathrm{int}(\mathcal U) \neq \varnothing$.
\item[\namedlabel{itm:5con}{(5)}] The control penalty $r \in C^1(\mathcal U)$ is strongly convex.
\item[\namedlabel{itm:growthcon}{(\textbf{\texttt{GC}})}] There exists a {$\boldsymbol{\pi}: \mathbb{X} \to \mathrm{int}(\mathcal U)$}, $\boldsymbol{\pi} \in L^\infty(\mathbb{X})$, a strongly convex $\mathcal V \in C^2(\mathbb{X})$ with bounded Hessian, and constants $C_1,C_2\in(0,\infty)$ such that
\begin{equation*}
(\boldsymbol{f}(\boldsymbol{x})+\boldsymbol{G}(\boldsymbol{x}) \boldsymbol{\pi}(\boldsymbol{x}))^\tr \nabla \mathcal V \leq C_1 - C_2 ( \stgcost(\boldsymbol{x}) + r(\boldsymbol{\pi}(\boldsymbol{x})) )
\end{equation*}
for all $\bx \in \mathbb{X}$.
\end{enumerate}

\end{assumption}
It is worth noting that \ref{itm:1con}--\ref{itm:5con} of Assumption~\ref{ass::blanket} are considered standard and typically do not impose restrictive limitations for practical applications. \ref{itm:growthcon} can be
viewed as a growth constraint imposed on $\boldsymbol{f},\boldsymbol{G},\stgcost,r$ as $\boldsymbol{x}$ tends to infinity. This criterion is often readily met in practical scenarios~\citep{houska2025convex}, as the constant $C_1$ can assume a
potentially large positive value.

\subsection{Fokker--Planck-Kolmogorov {PDE}s}

As long as Assumption~\ref{ass::blanket} holds, the martingale solution $\boldsymbol{X}_t$ in \eqref{eq::SDE} constitutes a time-homogeneous Markov diffusion process whose transition operator, $\Gamma_t(\boldsymbol{\pi}): L^2(\mathbb{X}) \to H^1(\mathbb{X})$, $t > 0$, maps the initial probability distribution $\rho_0 \in L^2(\mathbb{X})$ of $\boldsymbol{X}_0$ to the probability density function $\rho_t \in L^2(\mathbb{X})$ of $\boldsymbol{X}_t$,
\begin{equation*}
\textstyle \forall t \in (0,T], \qquad \rho_t = \Gamma_t(\boldsymbol{\pi}) \rho_0.
\end{equation*}
Note that $\Gamma_t(\boldsymbol{\pi})$ is for any given ${\boldsymbol{\pi}}\in \Pi$ and any $t > 0$ a bounded linear operator on $H^1(\mathbb{X})$~\citep{Hinze2009,Oksendal2000}. Its associated differential generator is given by
\begin{equation*}
\forall \rho \in C^\infty(\mathbb{X}), \qquad \mathcal L(\boldsymbol{\pi}) \rho\ \defeq \textstyle \lim_{t \to 0^+} \frac{\Gamma_t(\boldsymbol{\pi}) \rho - \rho}{t}.
\end{equation*}
This operator is also known as the Fokker--Planck-Kolmogorov (FPK) operator~\citep{Bogachev2015}. For differentiable feedback laws $\boldsymbol{\pi}$ and smooth test functions $\rho,V \in C^\infty(\mathbb{X})$, this operator is given by
\begin{equation}
\mathcal  L(\boldsymbol{\pi}) \rho = -\mathrm{div}(\rho(\boldsymbol{f}+\boldsymbol{G} \boldsymbol{\pi})) + \epsilon \Delta \rho
\end{equation}
For the case that $\boldsymbol{\pi} \in L^\infty(0,T;\Pi)$ is a potentially time-varying feedback law, the evolution of the probability density $\rho \in W(0,T;L^2(\mathbb{X}),H^1(\mathbb{X}))$ of the state of~\eqref{eq::SDE} satisfies --- if it exists --- the parabolic
FPK~\citep{Bogachev2015,Oksendal2000}
\begin{equation}
\label{eq::FPK}
\begin{array}
{rcl}
    \dot \rho_t &=& \mathcal L(\boldsymbol{\pi}_t) \rho_t = \mathcal{A} \rho_t + \mathcal{B} (\boldsymbol{\pi}_t \rho_t) \\
    \text{with} \quad \underset{t \to 0^+}{\lim}\rho_t &=& \rho_0
\end{array}
\end{equation}
in the weak sense, where $\rho_0 \in L^2(\mathbb X)$ denotes the probability density of $\bm{X}_0$. Note that $ \mathcal L(\boldsymbol{\pi})$ is affine in $\boldsymbol{\pi}$. Hence, we define the time-autonomous operators $\mathcal{A}$ and $\mathcal{B}$ by $\mathcal{A}\rho = -\mathrm{div} \left(\boldsymbol{f} \rho  \right) + \epsilon \Delta \rho $\, and $\, \mathcal{B}(\boldsymbol{\pi}\rho) = -\mathrm{div} \left(  \boldsymbol{G} \boldsymbol{\pi} \rho  \right)$ for all test functions $\rho \in C^\infty$.

If Assumption~\ref{ass::blanket} holds and the horizon is finite ($T<\infty$), then for any time-invariant feedback law $\boldsymbol{\pi} \in \Pi$ satisfying the last condition of Assumption~\ref{ass::blanket}, a unique weak solution to~\eqref{eq::FPK}
and a corresponding martingale solution $\boldsymbol{X}_t$ of~\eqref{eq::SDE} exist~\citep{houska2025convex,Oksendal2000}. In particular, the knowledge of a weak solution for at least one policy $\boldsymbol{\pi}$ suffices to guarantee
strict feasibility of~\eqref{eq::SOCP}~\citep{houska2025convex}.

\subsection{Convex Reformulation}\label{sec::ConvexReformulation}

By using the notion of weak solutions of parabolic FPKs, the original stochastic optimal control problem~\eqref{eq::SOCP} can be cast into a \emph{convex} PDE-constrained optimization
problem~\citep{Hinze2009,houska2025convex}, through a change of variables, $\boldsymbol{\nu} = \boldsymbol{\pi}{\rho}$,
\begin{eqnarray}
\label{eq::PDEOCP}
\mathcal J(T,\rho_0) & \defeq& \min_{\rho, \boldsymbol{\nu}} \int_0^T \int_{\mathbb{X}} \left(\stgcost+r\left(\frac{\boldsymbol{\nu}}{\rho} \right) \right) \rho \, \mathrm{d}\bx \, \mathrm{d}t \\
& &\text{s.t.}\ \left\{
\begin{array}
{ll}
\dot \rho  = \mathcal{A} \rho + \mathcal{B} \boldsymbol{\nu} & ~\text{on}~\mathbb{X}_T \\
\boldsymbol{\nu}\in \rho \mathcal U  & ~\text{on}~\mathbb{X}_T \\
\rho(0) = \rho_0 & ~\text{on}~\mathbb{X},
\end{array}\right. \notag
\end{eqnarray}
with $\mathbb{X}_T = (0,T) \times \mathbb{X}$. Since $\epsilon > 0$, the density $\rho_t = \Gamma_{t}(\boldsymbol{\pi}) \rho_0$ satisfies $\rho_t > 0$ for $t>0$, see \citet{Bogachev2015}.

\begin{proposition}[Houska 2025, Thm.~1]
\label{prop::PDEOCP}
If Assumption~\ref{ass::blanket} holds, $\ell$ has finite variance at the optimal steady state, and $\rho_0 \in L^2(\mathbb X)$ has bounded support, then~\eqref{eq::PDEOCP} admits a unique minimizer $(\rho^\star,\nu^\star)$ with
\begin{equation*}
\rho^\star \in W(0,T;L^2(\mathbb X),H^1(\mathbb{X})) \  \text{and} \ \nu^\star \in L^2(0,T;L^2(\mathbb{X}, \mathbb{R}^{n_u})).
\end{equation*}
Moreover, the optimal feedback law is given by $\boldsymbol{\pi}^\star = \frac{\boldsymbol{\nu}^\star}{\rho^\star}$ and $\rho^\star(t)$ corresponds to the probability density of the optimal state $\boldsymbol{X}_t^\star$.
\end{proposition}

\begin{remark}
If Assumption~\ref{ass::blanket} holds,
\begin{eqnarray}
\ell_\infty &=& \min_{\rho_\infty, \boldsymbol{\nu}_\infty} \int_{\mathbb{X}} \left(\stgcost+r\left(\frac{\boldsymbol{\nu}_\infty}{\rho_\infty} \right) \right) \rho_\infty \, \mathrm{d}\bx \, \notag \\
& &\text{s.t.}\ \left\{
\begin{array}
{ll}
0 = \mathcal{A} \rho_\infty + \mathcal{B} \boldsymbol{\nu}_\infty & ~\text{on}~\mathbb{X} \\
\boldsymbol{\nu}_\infty \in \rho_\infty \mathcal U  & ~\text{on}~\mathbb{X} \\
1 = \int_{\mathbb X} \rho_\infty \, \mathrm{d}\bx &
\end{array}\right. \notag
\end{eqnarray}
has a unique minimizer $(\rho_\infty^\star,\boldsymbol{\nu}_\infty^\star) \in H^1(\mathbb X) \times L^2(\mathbb X, \mathbb{R}^{n_u})$, see \citet[Cor.~1]{houska2025convex}. Moreover, under the additional assumption that the stage cost $\ell$ has a finite variance for the steady-state density, $\int (\ell-\ell_\infty)^2 \rho_\infty \, \mathrm{d}\bx < \infty$, we
have
\begin{equation*}
\ell_\infty =  \lim_{T \to \infty} \ \frac{1}{T} \ \mathcal J(T,\rho_0).
\end{equation*}
Because $\ell_\infty$ is constant, it can be interpreted as an ergodic limit that does not depend on the initial density $\rho_0$.
\end{remark}

\subsection{Hamilton--Jacobi-Bellman Equation}

Let us now connect \eqref{eq::PDEOCP} to the well-known stochastic Hamilton--Jacobi-Bellman (HJB) equation \citep{Fleming1993}.

Let $\mu$ be the ergodic measure at steady-state, such that $\mathrm{d} \mu = \rho_\infty dx$. Building on the strong duality based argument in \citet[Theorem 2]{houska2025convex}, we can construct a Riesz representation of the bounded
linear operator \begin{equation*}
 \mathcal J(T,\rho_0) = \int_{\mathbb X} V(0) \rho_0 \mathrm{d}\bx,
\end{equation*}
where $V(0)$ is the initial value of a strongly measurable function $V: [0,T] \to H_\mu^1(\mathbb{X})$, which turns out to be the weak solution of the HJB equation
\begin{equation}
\label{eq::HJB}
     -\dot V \! =  \!\mathcal A^{*} V \!+ \! \ell \! + \!\min_{\boldsymbol{\pi} \in \Pi}\left\{  r(\boldsymbol{\pi}) \!+\! \boldsymbol{\pi}^\tr \mathcal B^{*} V \right\}  ~\text{on}~ \ (0,T) \! \times \! \mathbb{X}
\end{equation}
with $V(T)= 0  $ on $\mathbb{X}$, thereby enabling the explicit computation of optimal controls and the corresponding value function for the stochastic control problem.

\subsection{Discretization of Time}

This section introduces a discrete-time approximation of the bounded linear operator $\mathcal J(T,\cdot)$. For this aim, we choose a step-size $h> 0$, a discrete-time horizon $H \in \mathbb{N}$ such that $T = Hh$, and define the discrete times $t_k = k h$ for $k \in \mathbb N$. By substituting a piecewise-constant approximation of the control law, $\boldsymbol{\pi}(t,\bx) = \boldsymbol{\pi}_k(\bx)$, for $t \in [t_k, t_{k+1}]$, we arrive at the discrete-time approximation
\begin{equation}
\label{eq:D-FPKOCP}
\begin{array}
{rcl}
\mathsf{J}_H(\rho_0) \ \defeq & \displaystyle\min_{\rho,\boldsymbol{\nu}}& \displaystyle\sum_{k=0}^{H-1} \int_{\IN} \left( \stgcost_{h} + r_{h}\left( \frac{{\boldsymbol{\nu}}_{\boldsymbol{k}}}{\rho_k} \right)  \right) \rho_k \, \mathrm{d}\bx \\
& \text{s.t.} & \left\{
\begin{array}
{rcl}
\rho_{k+1} &=& \mathsf{A}\, \rho_k +  \mathsf{B}\, {\boldsymbol{\nu}}_{\boldsymbol{k}} \quad\text{on}~\mathbb{X}   \\
{\boldsymbol{\nu}}_{\boldsymbol{k}} &\in& \rho_k \mathcal U,
\end{array}\right.
\end{array}
\end{equation}
which is defined for all $\rho_0 \in L^2(\IN)$. Here, $ \SFA$ and $ \SFB$ denote the Hankel transition operators of the linear differential equation~\eqref{eq::FPK},
\begin{equation}
\SFA \defeq \exp \left( \mathcal A {h} \right), \quad   \SFB \defeq \int_0^{h} \exp \left( \mathcal A \tau \right) \mathcal B \, \mathrm{d}\tau.
\end{equation}
Here, we define $\ell_h = h \ell$ and $r_h = h r$.
If Assumption~\ref{ass::blanket} holds, $\mathsf{J}_H$ is again a bounded linear operator on $L^2(\mathbb X)$ that admits a Riesz representation, $\mathsf{J}_H(\rho_0) = \int_\IN V_0 \rho_0 \, \mathrm{d}\bx$. Define the cost-to-go from step $k$ as $\mathsf{J}_k(\rho_k)$. Since each $\mathsf{J}_k$ is a bounded linear operator on $L^2(\IN)$, it admits a Riesz representation $\mathsf{J}_k(\rho_k) = \int_{\IN} V_k \, \rho_k \, \mathrm{d}\bx$ for some $V_k \in L^2(\IN)$, with terminal condition $V_H = 0$. To derive the backward recursion for $V_0$, we apply the Bellman optimality principle,
\begin{align*}
\mathsf{J}_k(\rho_k)
&= \min_{\boldsymbol{\nu}_k \in \rho_k \mathcal{U}}\int_{\IN}\left(
\stgcost_{h} + r_{h}\!\left(\frac{\boldsymbol{\nu}_k}{\rho_k}\right)
\right)\rho_k \, \mathrm{d}\bx
+\mathsf{J}_{k+1}(\SFA \rho_k + \SFB \boldsymbol{\nu}_k) \\
&= \min_{\boldsymbol{\nu}_k \in \rho_k \mathcal{U}}
    \int_{\IN}
    \left(\stgcost_{h} + r_{h}\!\left(\frac{\boldsymbol{\nu}_k}{\rho_k}\right)
    \right)\rho_k \, \mathrm{d}\bx +
    \int_{\IN} V_{k+1}\, (\SFA \rho_k + \SFB \boldsymbol{\nu}_k) \, \mathrm{d}\bx.
\end{align*}
where we substituted the Riesz representation of $\mathsf{J}_{k+1}$ and the dynamics $\rho_{k+1} = \SFA \rho_k + \SFB \boldsymbol{\nu}_k$ in the second line. Since $\boldsymbol{\nu}_k \in \rho_k \mathcal U$, there exists a measurable policy
$\boldsymbol{\pi}_k \in \mathcal U$ such that $\boldsymbol{\nu}_k = \rho_k \boldsymbol{\pi}_k$. Then, we combine all terms under a single integral to obtain
\begin{equation*}
\mathsf{J}_k(\rho_k)
= \min_{\boldsymbol{\pi}_k \in \Pi}\int_{\IN}
\Bigl(
\stgcost_{h} + r_{h}(\boldsymbol{\pi}_k) + \SFA^{*} V_{k+1}
+
\boldsymbol{\pi}_k^\tr \SFB^{*} V_{k+1}
\Bigr)\rho_k
\, \mathrm{d}\bx.
\end{equation*}
where we applied the adjoint relation
$
\int_{\IN} V_{k+1}\, \SFA \rho_k \, \mathrm{d}\bx
= \int_{\IN} (\SFA^{*} V_{k+1})\, \rho_k \, \mathrm{d}\bx
$ and analogously for $\SFB$.
Identifying the integrand with $V_k$ and dropping the time index on the optimized variable $\boldsymbol{\pi}$ yields the following operator dynamic programming recursion, from which $V_0$ is obtained. In detail, we have
\begin{equation}
\label{eq:OTHJB}
\begin{aligned}
V_k &= \SFA^{*} V_{k+1} + \ell_h
     + \min_{\boldsymbol{\pi}\in\Pi}\!\left\{ r_h(\boldsymbol{\pi}) + \boldsymbol{\pi}^{\tr}\SFB^{*} V_{k+1} \right\}
     && \text{on } \mathbb{X},\\
V_H &= 0 && \text{on } \mathbb{X},
\end{aligned}
\end{equation}
which can be interpreted as a time-discretization of \eqref{eq::HJB}.

\begin{remark}\label{rem:timeInvariantPolicy}
For finite-horizon optimal control problems, the above time-discretization leads to a discretization error, since the optimal density and control law are, in general, not time-invariant. Nevertheless, in the ergodic limit, we still have
\begin{equation*}
\ell_\infty = \lim_{T \to \infty} \ \frac{1}{T} \ \mathcal J(T,\rho_0) = \lim_{H \to \infty} \ \frac{1}{hH} \ \mathsf{J}_H(\rho_0)
\end{equation*}
for any $h > 0$ and any $\rho_0 \in L^2(\mathbb X)$, because the ergodic feedback law, given by $\boldsymbol{\pi}_\infty = \frac{\boldsymbol{\nu}_{\boldsymbol{\infty}}}{\rho_\infty}$, is time-invariant.
\end{remark}

\section{Optimal Control Meets RKHS Embeddings}
\label{sec:optctrlKME}
The objective of this section is to propose a data-driven approach for approximating the solution of the stochastic optimal control problem.

\subsection{Reproducing Kernel Hilbert Spaces}

Reproducing Kernel Hilbert Spaces (RKHS) are widely accepted as the key to nonparametric statistical learning in function spaces \citep{IngoSteinwart2008SupportMachines}. Let $k: \mathbb{X} \times \mathbb{X} \to \mathbb R$ be a non-negative, symmetric, and
continuous kernel and let $(\spIN,\innerprod{\cdot,\cdot}_{\spIN})$ be an associated RKHS. This means that $(\spIN,\innerprod{\cdot,\cdot}_{\spIN})$ is a Hilbert space, whose inner product satisfies the reproducing property $\psi(\bx) = \innerprod{k(\cdot,\bx), \psi}_{\spIN}$ for all $\bx \in \IN$ and all observables
$\psi \in \spIN$. In the following, we assume that $\spIN$
is a dense subspace of $H^1(\IN)$, $k$ is universal,\footnote{This is a property satisfied by many common kernels such as Gaussian and Mat\'ern, see \citet{IngoSteinwart2008SupportMachines}.}  and $k(\bx, \bx') < \infty$
for all $\bx,\bx' \in \mathbb X$.

\subsection{{RKHS} Embeddings}

We introduce the following embedding of a probability density $\rho \in L^2(\IN)$ into an RKHS $\spIN$ as
\begin{equation}
\KME{\SFX}: L^2(\IN) \to \spIN, \quad \rho \to \int_{\IN} k(\cdot, \bx) \rho(\bx)\,\mathrm{d}\bx
    \in \spIN \label{eq:kme_p}
\end{equation}
It can be interpreted as the expected value of the feature map, namely $ \mathcal{E}_{\SFX} \rho = \mathsf{E}_{\boldsymbol{X} \sim \rho}[k(\cdot, \boldsymbol{X})]$. Equivalently, $ \mathcal{E}_{\SFX}$ can be viewed as the kernel mean embedding (KME) \citep{Muandet2017} of the probability measure $\mu$
with density $\rho$ with respect to the Lebesgue measure, $d\mu(\bx) = \rho(\bx) \mathrm{d}\bx$. We also define the embedding of the controls $\boldsymbol{\nu} \in L^2 (\IN,  \mathbb{R}^{n_u})$ into an RKHS $\mathbb{R}^{n_u}\!\otimes\!\spIN$ by
\begin{equation}
\KME{\SFX, \SFU} \, \boldsymbol{\nu} \defeq \int_{\IN} \boldsymbol{\nu}(\bx)\!\otimes\!k(\cdot, \bx)\, \mathrm{d} \bx
    \in \mathbb{R}^{n_u}\!\otimes\!\spIN. \label{eq:kme_u}
\end{equation}

Next, we embed and approximate the dynamics in \eqref{eq:D-FPKOCP} as
\begin{equation}
\mathcal{E}_{\SFX} \rho_{k+1}= \mathcal{E}_{\SFX} (\SFA \rho_k + \SFB \boldsymbol{\nu}_k)  \approx A \mathcal{E}_{\SFX} \rho_k + B\KME{\SFX, \SFU} \boldsymbol{\nu}_k \label{eq:KME_dyn},
\end{equation}
where the left-hand equality uses the linearity of the embedding~\eqref{eq:kme_p}. The approximation assumes the existence of Hilbert--Schmidt operators $A \in \mathrm{HS}(\spIN, \spIN)$ and $B \in \mathrm{HS}(\mathbb{R}^{n_u}\!\otimes\!\spIN, \spIN)$ such that $\KME{\SFX} \SFA \rho \approx A \KME{\SFX} \rho$ and
$\KME{\SFX} \SFB \boldsymbol{\nu} \approx B \KME{\SFX, \SFU} \boldsymbol{\nu}$, for all $\rho \in L^2(\IN)$ and $\boldsymbol{\nu} \in L^2(\IN, \mathbb{R}^{n_u})$. That is, the transition of embedded densities under the transition operators $(\SFA, \SFB)$ is approximated by the operators $(A, B)$ acting directly on the RKHS
embeddings. This is well-motivated because a universal kernel ensures that the feature space $\spIN$ is rich enough to distinguish any two distinct densities through their embeddings (i.e.,~the embeddings
\eqref{eq:kme_p} and \eqref{eq:kme_u} are injective~\citep{Muandet2017}), so the embedded representation retains all the relevant information about the underlying dynamics. The operators $(A, B)$ are then learned from data to
minimize this approximation error in the next paragraph.

\subsection{Data-Driven Operator Learning}

In practice, we do not have access to population (infinite data) levels to work directly with the embeddings introduced previously. However, we assume the existence of (finite) training samples of the form
\begin{equation}
\label{eq:data}
    \Set{D}_{n}= \bigl\{(\boldsymbol{x}^{(i)},\boldsymbol{u}^{(i)}), \, \boldsymbol{x}_{+}^{(i)}, \, \stgcost(\boldsymbol{x}^{(i)}) \bigr\}_{i{\in}[n]},
\end{equation}
where the states and controls are selected according to a strictly positive probability measure $\mu \in M_+(\IN \times \mathbb{U})$. Although in practice the training samples are collected on a compact subset of $\IN \times \mathbb U$, we retain this notation for simplicity. Here, $\boldsymbol{x}_{+}^{(i)}$ denotes the successor state of $\boldsymbol{x}^{(i)}$ under the sampling time $h > 0$, corresponding to the discrete-time observations $(\boldsymbol{x}_k, \boldsymbol{u}_k, \boldsymbol{x}_{k+1})$ of the SDE~\eqref{eq::SDE} introduced in Section~\ref{sec:informalProblem}. Note that the stage cost values $\stgcost(\boldsymbol{x}^{(i)})$ may either be sampled (as shown) or computed directly if $\stgcost$ is known a priori.

We define the dataset-dependent empirical embeddings $\EKME{\SFX},\EKME{+}:\mathbb{R}^n \to \spIN$ and $\EKME{\SFX, \SFU}:  \mathbb{R}^{n}\to \mathbb{R}^{n_u} \otimes \spIN$ by
\begin{align}
\EKME{\SFX} \boldsymbol{w} &\defeq \textstyle\frac{1}{\sqrt{n}}\lilsum_{i=1}^n w_i k(\cdot, \stmbxi), \label{eq:ekmeX} \\
    \EKME{+} \boldsymbol{w} &\defeq \textstyle\frac{1}{\sqrt{n}}\lilsum_{i=1}^n w_i k(\cdot, \boldsymbol{x}_{+}^{(i)}), \label{eq:ekmeXplus} \\
    \EKME{\SFX, \SFU} \boldsymbol{w} &\defeq \textstyle\frac{1}{\sqrt{n}}\lilsum_{i=1}^n w_i\bui\!\otimes\!k(\cdot, \stmbxi) \label{eq:ekmeU}
\end{align}
for all $\boldsymbol{w} \in \mathbb{R}^n$. Let $(\boldsymbol{e}_i)_{i=1}^n$ denote the canonical basis of $\mathbb{R}^n$. By substituting the empirical embeddings \eqref{eq:ekmeX}--\eqref{eq:ekmeU} in \eqref{eq:KME_dyn}, we can write the empirical residual as
\begin{align}
\!\!\! \hat{\mathsf{r}}_i &=  \left(\EKME{+}  - ( A \EKME{X} + B \EKME{X, U} )\right)\boldsymbol{e}_i \\
   &=  \textstyle\frac{1}{\sqrt{n}}\left(k(\cdot, \boldsymbol{x}_{+}^{(i)}) {-} \left( A k(\cdot, \stmbxi) {+} B \left(\bui\!\otimes\! k(\cdot, \stmbxi)\right)\right)\right) \\
   &= \textstyle\frac{1}{\sqrt{n}}k(\cdot, \boldsymbol{x}_{+}^{(i)}) - \textstyle\frac{1}{\sqrt{n}}P\left( \concat{1}{\boldsymbol{u}^{(i)}} \otimes k(\cdot, \stmbxi)\right),
   \label{eq:residual}
\end{align}
where we defined the operator $P$ $\in \mathrm{HS}(\spOUT, \spIN)$, with the RKHS $\spOUT = \spIN \oplus(\mathbb{R}^{n_u} \otimes \spIN)$. For brevity, we also introduce the shorthand $\boldsymbol{z} \defeq \concat{\boldsymbol{x}}{\boldsymbol{u}}$ so that $\Set{Z}\defeq \Set{X}{\times}\Set{U}$.

Given the training samples \eqref{eq:data}, we can obtain an estimator \citep{Muandet2017} by minimizing the squared RKHS norm error over the residuals $\hat{\mathsf{r}}_i$ of \eqref{eq:residual}, thus solving
\begin{equation}
\label{eq:Estimator1}
\widehat{P}  = \argmin_{P \in \mathrm{HS}(\spOUT, \spIN)} {\textstyle\frac{1}{n}} \sum_{i=1}^n\big\Vert k(\cdot, \bx_+^{(i)}) - P\left( \concat{1}{\boldsymbol{u}^{(i)}} \otimes k(\cdot, \bx^{(i)})\right) \big\Vert ^2_{\spIN} + \gamma \|  P \| ^2_{\mathrm{HS}},
\end{equation}
where $\gamma$ is a regularization, which ensures numerical stability and prevents overfitting. In addition to \eqref{eq:ekmeX}--\eqref{eq:ekmeU}, we introduce the sampling operators
\begin{align}
\widehat{\mathcal{E}}_{\SFZ}^{*}&:\spOUT\to\Set{R}^n,\quad \widehat{\mathcal{E}}_{\SFZ}^{*}k_{\mathsf{Z}}(\cdot, \bz) = \textstyle\frac{1}{\sqrt{n}}\left[{k}_{\SFZ}(\bz, \bz^{(i)})\right]_{i \in [n]}, \\
    \widehat{\mathcal{E}}_+^{*}&:\spIN\to\Set{R}^n, \quad \widehat{\mathcal{E}}_+^{*} k(\cdot, \bx) =  \textstyle\frac{1}{\sqrt{n}}\left[k(\bx, \bx_+^{(i)})\right]_{i \in [n]},
\end{align}
where ${k}_{\SFZ}(\boldsymbol{z},\boldsymbol{z}^\prime) ={k}(\boldsymbol{x}^\prime,\boldsymbol{x})(1{+}\innerprod{\boldsymbol{u}, \boldsymbol{u}^\prime})$, see \citet{bevanda2025nonparametric}. Note that the above sampling operators are the adjoints of the corresponding empirical embeddings $ \widehat{\mathcal{E}}_{\SFZ}$ and $ \widehat{\mathcal{E}}_+$. Applying the kernel trick, we can
obtain the Gram matrices $(\widehat{\mathcal{E}}_{\SFX}^{*} \widehat{\mathcal{E}}_{\SFX})_{ij} \defeq (\Kx)_{ij}=\textstyle \frac{1}{n} k(\bx^{(i)},\bx^{(j)})$ and $(\widehat{\mathcal{E}}_{\SFZ}^{*} \widehat{\mathcal{E}}_{\SFZ})_{ij} \defeq (\boldsymbol{K})_{ij}= \textstyle \frac{1}{n} {k}_{\SFZ}(\bz^{(i)},\bz^{(j)})$. The regularized inverse of this Gram matrix is denoted by $\boldsymbol{K}^{-1}_{\gamma} = (\widehat{\mathcal{E}}_{\SFZ}^{*} \widehat{\mathcal{E}}_{\SFZ} + \gamma \boldsymbol{I})^{-1} = (\boldsymbol{K}  + \gamma \boldsymbol{I})^{-1} $.
To solve~\eqref{eq:Estimator1}, we note that the sum of squared residuals associated with \eqref{eq:residual} equals a single Hilbert--Schmidt norm. Indeed, $\sum_{i=1}^n \| \hat{\mathsf{r}}_i \| _\spIN^2 =  \sum_{i=1}^n \|  (\widehat{\mathcal{E}}_+ - P\widehat{\mathcal{E}}_{\SFZ} ) \boldsymbol{e}_i  \| _\spIN^2 =\| \widehat{\mathcal{E}}_+ - P\widehat{\mathcal{E}}_{\SFZ}\| ^2_{\mathrm{HS}}$, by orthonormality of $\{\boldsymbol{e}_i\}_{i\in[n]}$ and the definition of the Hilbert--Schmidt norm. Hence, the first-order optimality condition of this resulting Kernel Ridge Regression (KRR) leads to
\begin{equation*}
\widehat{P}\bigl(\widehat{\mathcal{E}}_{\SFZ} \widehat{\mathcal{E}}_{\SFZ}^{*} + \gamma \mathsf{I}\bigr) = \widehat{\mathcal{E}}_+ \widehat{\mathcal{E}}_{\SFZ}^{*},
\end{equation*}
so that $\widehat{P} = \widehat{\mathcal{E}}_+ \widehat{\mathcal{E}}_{\SFZ}^{*} (\widehat{\mathcal{E}}_{\SFZ} \widehat{\mathcal{E}}_{\SFZ}^{*} + \gamma \mathsf{I})^{-1}$. Applying the push-through identity gives $ \widehat{P} = \widehat{\mathcal{E}}_+(\widehat{\mathcal{E}}_{\SFZ}^{*} \widehat{\mathcal{E}}_{\SFZ} + \gamma \boldsymbol{I})^{-1} \widehat{\mathcal{E}}^{*}_{\SFZ} = \widehat{\mathcal{E}}_+\boldsymbol{K}^{-1}_{\gamma} \widehat{\mathcal{E}}^{*}_{\SFZ}$. Finally, under the decomposition $\spOUT = \spIN \oplus (\mathbb{R}^{n_u} \otimes \spIN)$, the adjoint $\widehat{\mathcal{E}}_{\SFZ}^{*} : \spOUT \to \mathbb{R}^n$ decomposes accordingly, with first component $\widehat{\mathcal{E}}_{\SFX}^{*}$ and second component $\widehat{\mathcal{E}}_{\SFX,\SFU}^{*}$, yielding
\begin{equation}
\widehat{P} = [{\widehat{A}}~~{\widehat{B}}]=  [ \widehat{\mathcal{E}}_+  \boldsymbol{K}^{-1}_{\gamma} \widehat{\mathcal{E}}_{\SFX}^{*} ~~~\widehat{\mathcal{E}}_+ \boldsymbol{K}^{-1}_{\gamma} \widehat{\mathcal{E}}_{\SFX,\SFU}^{*}]\label{eq:estimators},
\end{equation}
which is dual to learning control Koopman Operators (cKOR) in \citet{bevanda2025nonparametric}.

The regression formulation above resembles the Extended Dynamic Mode Decomposition (EDMD) framework \citep{williams2015data}, but uses feature maps in an RKHS instead of an explicit observable dictionary. This obviates the
need for explicit tensor products between the dictionary and controls \citep{STRASSER20232257}. Furthermore, with universal kernels, the resulting kernel methods allow arbitrarily accurate learning of these operators.

\subsection{Data-Driven Optimal Control}

In this section, we aim to develop a data-driven approximation of \eqref{eq:OTHJB}, using the estimators $\widehat{A}, \widehat{B}$ from the previous section. Recall that the Riesz representer $V_0$ is obtained as the (weak) solution of the discretized stochastic HJB equation in \eqref{eq:OTHJB}. Let us define the Fenchel conjugate
\begin{align}
D_r( \boldsymbol{\lambda} ) &\defeq \min_{\boldsymbol{u} \in \mathcal U}\left\{  r_h(\boldsymbol{u}) + \boldsymbol{\lambda}^\tr \boldsymbol{u} \right\} \label{eq:dualfct} \\
    \text{with} \quad \boldsymbol{u}^\star( \boldsymbol{\lambda}) &\defeq \underset{\boldsymbol{u} \in \mathcal U}{\argmin}\left\{ r_h(\boldsymbol{u}) + \boldsymbol{\lambda}^\tr \boldsymbol{u} \right\}, \label{eq:minimizer}
\end{align}
where, as in \eqref{eq:OTHJB}, we set  $r_h = h r$. Since we assume that $r$ is strongly convex, while $\mathcal U$ is closed, convex, and has a non-empty interior, the parametric minimizer $\boldsymbol{u}^\star$ is unique and Lipschitz continuous. To derive a surrogate of \eqref{eq:OTHJB} in the RKHS $\spIN$, we provide the cost $\ell_h = h \ell$ and the Fenchel conjugate ${D}_r$ representations via $\EKME{\SFX}$ \eqref{eq:ekmeX}
\begin{equation*}
        \widehat{\ell}_{h}  \defeq  \widehat{\mathcal{E}}_{\SFX} \boldsymbol{\stgcost}_{h}, \quad \widehat{D}_r(\cdot ) \defeq \widehat{\mathcal{E}}_{\SFX}\boldsymbol{D}_r (\cdot),
\end{equation*}
respectively, where $\SFD_r(\cdot): \mathbb{R}^{n n_u} \to \mathbb{R}^n$ and $\boldsymbol{\stgcost}_{h} \in \mathbb{R}^n$. Substituting the adjoints of the learned operator model $(\widehat{A}, \widehat{B})$ from \eqref{eq:estimators} into \eqref{eq:OTHJB} then yields the following RKHS embedded dynamic programming recursion for $\EV_k \in \spIN$
\begin{equation}
\EV_k =  \widehat{A}^{*}  \EV_{k+1} + \widehat{\ell}_{h} + \widehat{D}_r( \widehat{B}^{*}   \EV_{k+1} ), \quad \EV_H = 0,\label{eq:HJBrecursiononH}
\end{equation}
where $h > 0$ and $H$ denote the discrete-time step and horizon, respectively.
While \eqref{eq:HJBrecursiononH} is intractable in this form, under the finite-data setting \eqref{eq:data} we obtain finite-rank $\widehat{A}$ and $\widehat{B}$ in \eqref{eq:estimators}. Next, we insert their adjoints $( \widehat{A}^{*},  \widehat{B}^{*})$ into \eqref{eq:HJBrecursiononH} and test $\EV_{k+1}$ against $k(\cdot, \bx)$, which yields
\begin{equation*}
\innerprod{\widehat{A} \featx{\bx}, \EV_{k+1}}_\spIN \! \! =  \!\innerprod{\widehat{\mathcal{E}}_{\SFX}^{*}\featx{\bx}, \boldsymbol{K}^{-1}_{\gamma} \!\boldsymbol{K}_{+} \SFV_{k+1}}_{\mathbb{R}^n} \! \!= \! \innerprod{\boldsymbol{k}(\bx), \boldsymbol{A} \SFV_{k+1}}_{\mathbb{R}^n}
\end{equation*}
with $(\boldsymbol{K}_{+})_{ij} = (\widehat{\mathcal{E}}_+^{*} \widehat{\mathcal{E}}_{\SFX})_{ij} = \textstyle\frac{1}{n} k(\boldsymbol{x}_+^{(i)},\bx^{(j)})$ and
\begin{equation*}
    \innerprod{\widehat{B} (\bu \otimes \featx{\bx}), \EV_{k+1}}_\spIN = \innerprod{\mathrm{diag}(\boldsymbol{U}\bu)\boldsymbol{k}(\bx), \boldsymbol{A} \SFV_{k+1}}_{\mathbb{R}^n},
\end{equation*}
where we used $\widehat{\mathcal{E}}_{\SFX, \SFU}^{*}(\bu \otimes \featx{\bx}) = \mathrm{diag}(\boldsymbol{U}\bu) \boldsymbol{k}(\bx)$ and defined $\boldsymbol{U} \coloneq  [\bu^{(1)},\dots,\bu^{(n)}]^\top$. Thus, we obtain \emph{tractable} matrix representations
\begin{equation}
\boldsymbol{A} = \boldsymbol{K}^{-1}_{\gamma} \boldsymbol{K}_{+}, \quad \boldsymbol{B} =[\mathrm{diag}(\boldsymbol{U} \boldsymbol{e}_{i})\boldsymbol{K}^{-1}_{\gamma}\boldsymbol{K}_{+}]_{i \in [n_u]}.
    \label{eq:MatrixAandB}
\end{equation}
Additionally, we define the RKHS representations
\begin{equation}
\widehat{\ell}_{h}  \defeq  \widehat{\mathcal{E}}_{\SFX} \boldsymbol{\stgcost}_{h}, \quad \widehat{D}_r( \widehat{B}^{*}   \EV_{k+1} ) \defeq   \widehat{\mathcal{E}}_{\SFX} \SFD_r( \boldsymbol{B} \SFV_{k+1} ) ,\label{eq:FenchelAndCostinH}
\end{equation}
where  $\SFD_r( \boldsymbol{B} \SFV_{k+1} ) = \boldsymbol{K}^{-1}_{\gamma} [D_r([\widehat{B}^{*} \EV_{k+1}](\stmbxi))]_{i \in [n]}$ and $\boldsymbol{\stgcost}_{h} = \boldsymbol{K}^{-1}_{\gamma} [\stgcost_h^{(i)}]_{i \in [n]}$ are the Kernel Ridge Regression (KRR) solutions corresponding to the evaluated dual terms on the dataset \eqref{eq:data} and the sampled stage costs, respectively. Since elements of $\mathbb{R}^{n_u} \otimes \spIN$ are equivalent to $\spIN^{n_u}$, $[\widehat{B}^{*} \EV_{k+1}](\stmbxi) =\innerprod{\boldsymbol{I}_{n_u} \otimes \boldsymbol{k}(\stmbxi), \,  \boldsymbol{B} \SFV_{k+1}}_{\mathbb{R}^{n n_u}}$.
These derivations lead to an exact and tractable recursion with matrix representations, as detailed in the following result.
\begin{theorem}
Let  $\widehat{A}, \widehat{B}$ be the estimators in \eqref{eq:estimators}, and define
$ \widehat{\ell}_{h}$ and $\widehat{D}_r( \widehat{B}^{*}   \EV_{k+1} )$ as in \eqref{eq:FenchelAndCostinH}. Then the operator recursion \eqref{eq:HJBrecursiononH} evolves in a finite-$n$-dimensional subspace, and is equivalent to the Kernel HJB (KHJB) recursion
\begin{equation}
\begin{aligned}
\label{eq:KHJB-system}
       &\SFV_{k} = \boldsymbol{A} \SFV_{k+1}  + \boldsymbol{\stgcost}_{h}   + \SFD_r( \boldsymbol{B} \SFV_{k+1} )  \\ \; \qquad \text{s.t.} \;   & \SFV_H \!=\ \mathbf{0}
\end{aligned}
\end{equation}
with $\SFV_k \in \mathbb{R}^n$ and matrices $(\boldsymbol{A}, \boldsymbol{B})$ given in \eqref{eq:MatrixAandB}.
\end{theorem}

\begin{proof}
This statement follows from the finite-rank representation for $\widehat{A}, \widehat{B}$. The first step to obtain \eqref{eq:KHJB-system} consists in writing the empirical embedding $\EV_k = \widehat{\mathcal{E}}_{\SFX} \SFV_{k}$ and leverage the reproducing property $\EV_k(\bx) = \innerprod{\featx{\bx}, \EV_k} = \innerprod{\widehat{\mathcal{E}}_{\SFX}^{*} \featx{\bx},\SFV_k} = \innerprod{\boldsymbol k(\bx),\SFV_k }$ with the sampled canonical feature map $\boldsymbol k(\bx) = \textstyle \frac{1}{\sqrt{n}} [k(\bx,  \stmbxi) ]_{i \in [n]}$.
Then \eqref{eq:KHJB-system} follows directly from substituting $(\widehat{A}^{*}, \widehat{B}^{*})$ along with \eqref{eq:FenchelAndCostinH} in \eqref{eq:HJBrecursiononH}.
\end{proof}

We call the recursion \eqref{eq:KHJB-system} starting at $\SFV_H$ the discrete-time KHJB, because it can be interpreted as a discretized version of the original HJB \eqref{eq::HJB}, for $T=Hh$. It can be used to approximate the (true) optimal value function $V^\star$ --- that is, the optimal co-state of \eqref{eq::PDEOCP} --- via the reproducing property
$\EV^\star(\bx) = \innerprod{\boldsymbol{k}(\bx), \SFV^\star}$ with $\SFV^\star = \SFV_0$.
This approximation introduces two main sources of error: the \emph{discretization error} and the \emph{learning error} incurred during the estimation of the operator model in \eqref{eq:estimators}. Note that for the latter, under standard regularity assumptions and using universal kernels, the data-driven estimators of the operators in \eqref{eq:estimators} are consistent and converge to the true operators
in the data limit; see \citet{Kostic2023KoopmanEigenvalues}. A more in-depth error analysis, however, is beyond the scope of this work.

To return to our problem formulation and the objective of finding an optimal feedback control law $\boldsymbol{\pi}^\star$ minimizing the ergodic cost in the Stochastic Optimal Control Problem~\eqref{eq::SOCP}, an approximate solution can be readily obtained as
\begin{equation}
\widehat{\boldsymbol{\pi}}^\star(\bx)
    = \boldsymbol{u}^\star\!\bigl(\innerprod{\boldsymbol{I}_{n_u} \otimes \boldsymbol{k}(\bx),\, \boldsymbol{B}\SFV^\star}_{\mathbb{R}^{n n_u}} \bigr).
\end{equation}

\begin{algorithm}[!ht]
\caption{Nonparametric Optimal Control}
\label{alg:ddoc}
\begin{algorithmic}[0]
\Require Data $\Set{D}_{n}$ \eqref{eq:data}, minimizer $\bm{u}^\star(\cdot)$ \eqref{eq:minimizer}, kernel $k$, horizon $H$ and sampling time ${h} > 0$; regularization parameter $\gamma > 0$, tolerance $\mathsf{tol}>0$.
\smallskip

\State Set $\bm{U} := [\bu^{(1)},\dots,\bu^{(n)}]^\top$
\State Set $\bm{k}(\bx) = \textstyle\frac{1}{\sqrt{n}} [k(\bx,  \bxi) ]_{i \in [n]}$ \Comment{$\bm{k}(\bx) \in \mathbb{R}^{n}$}
\smallskip \smallskip \smallskip
\State \texttt{/* NONPARAMETRIC MODEL LEARNING */}
\smallskip
\State Compute $(\Kx)_{ij} \defeq \textstyle\frac{1}{n} k(\bx^{(i)},\bx^{(j)})$
\State Compute $(\bm{K}_+)_{ij}\defeq \textstyle\frac{1}{n} k(\bm{x}_+^{(i)},\bx^{(j)})$
\State Compute $\bm{K}^{-1}_{\gamma} \defeq \left(\Kx+\Kx\odot \bm{U} \bm{U}^\tr + \gamma \bm{I}\right)^{-1} $

\State Compute $ \bm{A}  = \bm{K}^{-1}_{\gamma} \bm{K}_+ $
\State Compute $\bm{B} = [\mathrm{diag}(\bm{U} \bm{e}_{i}) \bm{A}]_{i \in [n_u]}$
\smallskip \smallskip \smallskip
\State \texttt{/* OPTIMAL CONTROL (KHJB) */}
\smallskip
\State Set $D_r(\bm{\lambda}) = r(\bm{u}^\star(\bm{\lambda}))h  + \bm{\lambda}^\tr \bm{u}^\star(\bm{\lambda}) $
\State $\bm{\lambda}(\bxi) = \innerprod{\bm{I}_{n_u} \otimes \bm{k}(\bxi), \,  \bm{\lambda}}\in \mathbb{R}^{n_u}$ \Comment{$\bm{\lambda} \in \mathbb{R}^{n n_u}$}
\State Set $\bm{D}_r(\bm{\lambda})  = \bm{K}^{-1}_{\gamma} [D_r(\bm{\lambda}(\bxi))]_{i \in [n]}$
\State Compute $\bm{\stgcost}_{h} = h \, \bm{K}^{-1}_{\gamma} [\stgcost(\bm{x}^{(i)})]_{i \in [n]} $
\smallskip
\State  Set $\SFV_H \gets\bm{0}, \;  k \gets H$
\Repeat
    \State $\SFV_{k-1} \ =\bm{A} \SFV_{k} + \bm{\stgcost}_h   + \SFD_r( \bm{B} \SFV_{k} )$
    \State $k\gets k -1$
\Until $k=0$ or $ \| \SFV_{k} -\SFV_{k+1} \|\leq  \mathsf{tol}$
\State $\SFV^\star =\SFV_{k}  $
\Output
\State  \!\!\!\!\!\!$\EV^{\star}(\bx) {=} \bm k(\bx)^\tr \SFV^\star$ and $\widehat{\bm{\pi}}^\star(\bx) {=} \bm{u}^\star( \innerprod{\bm{I}_{n_u} \otimes \bm k(\bx), \bm{B}\SFV^\star })$
\EndOutput
\end{algorithmic}
\end{algorithm}

\begin{table}[!htbp]\centering
\caption{Computational complexity of Algorithm~\ref{alg:ddoc}. {Here, $n$ denotes the number of training samples, $n_x$ the state dimension, $n_u$ the control dimension, and $H$ the discrete-time horizon.}}\label{tab:complexity}

\begin{tabular}{c|c}
\toprule
  Task
        &  \text{Complexity} \\
        \midrule
           Regression
           & $\order{ n^3 + n^2(n_x + n_u)}$ \\
           Control (\texttt{KHJB})
           & $\order{H (n^2  n_u +  n n_u^2)}$
\\
\bottomrule
\end{tabular}
    \end{table}

\section{Implementation and Numerical Examples}\label{sec:Numerical}

\subsection{The {KHJB} Algorithm}
The implementation follows Algorithm \ref{alg:ddoc}, with its computational complexity summarized in Table~\ref{tab:complexity}.
The main bottleneck is the inversion of the $n\times n$ matrix defining $\boldsymbol{K}^{-1}_{\gamma}$, which scales as $\mathcal{O}(n^3)$.
Each control-update step requires evaluating the Fenchel conjugate of the control penalty ($\mathcal{O}(n_u^2)$) and computing $\boldsymbol{B}\boldsymbol{v}$ ($\mathcal{O}(n^2 n_u)$). It is also worth noting that the complexity scales \emph{linearly} with
the state dimensionality $n_x$.  In large-scale settings, this dependence can be mitigated by low-rank kernel approximations and sketching, such as Nystr\"om methods and sketched operator regression
\citep{Caldarelli2024,Rudi2017,meanti2023estimating,bevanda2025nonparametric}, which reduce the dominant regression cost to $\mathcal O(nm(n_x+n_u)+nm^2+m^3)$ for landmark data points $m\ll n$, with the three terms corresponding to kernel
evaluation, reduced matrix assembly, and reduced-system inversion, respectively. The nonparametric operator regression viewpoint also remains compatible with additional tools for optimizing computational complexity,
including principal component or reduced rank regression surrogates in RKHS \citep{Kostic2022LearningSpaces}. These ideas also point to online implementations with a fixed computational budget. In that setting, the
operator estimator and the KHJB recursion can be updated incrementally while keeping $\boldsymbol{A}$ and $\boldsymbol{B}$ at constant size as new data arrive. This is conceptually similar to recursive estimation and control
pipelines familiar from Linear Quadratic Gaussian (LQG) settings, but now for nonlinear operator models.

\subsection{Numerical Examples}

In this section, we present numerical examples implemented in \textsf{Julia} to evaluate the accuracy of the data-driven generated global optimal control law using Algorithm \ref{alg:ddoc}.

Our numerical examples are based on Radial Basis Function (RBF) kernels of the form
$k(\boldsymbol{x},\boldsymbol{y}) = \exp{(\textstyle -\nicefrac{\| \boldsymbol{x} -\boldsymbol{y}\| ^ 2}{\sigma^2})}$. The RBF scale $\sigma$ is tuned separately for each benchmark, using a grid search based on the prediction performance. For the 1D and 2D examples, we set  $\gamma = 10^ {-8}$ and $\epsilon = 0.02$ and collect data using a uniform random scalar input $u$ on the domain $\mathbb U = [-1, 1]$.

\begin{figure*}[ht]
    \centering
    \setlength{\tabcolsep}{2pt}
    \begin{tabular}{cc}
        \subfloat[\texttt{S1}, $\pi^{\star}(x) = -\sqrt{2}\,x$]{
            \includegraphics[width=0.44\textwidth]{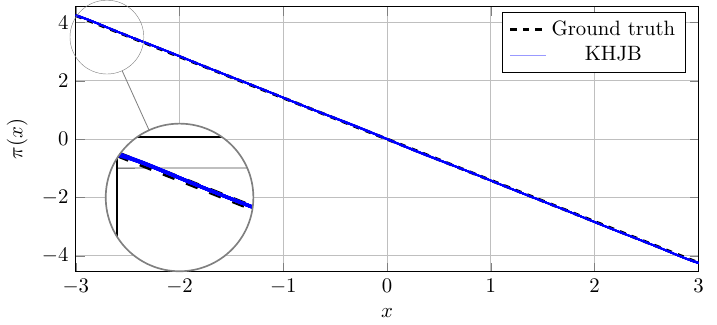}
            \label{fig:subfig1}
        } &
        \subfloat[\texttt{S2}, $\pi^{\star}(x) = -\ln{(x^2)}\,x$]{
            \includegraphics[width=0.42\textwidth]{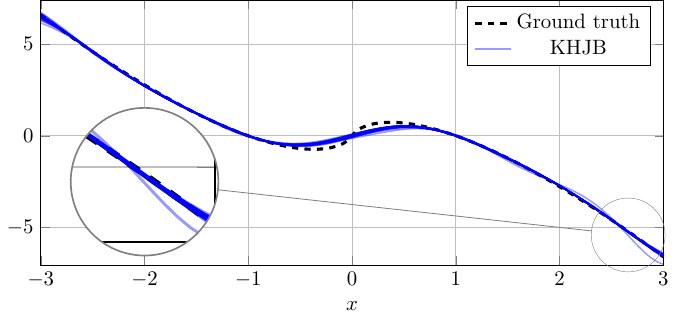}
            \label{fig:subfig2}
        } \\[-0.3em]
        \subfloat[\texttt{S3}, $\pi^{\star}(x) = -\tfrac{1}{2}x - \sin(2x)\,x$]{
            \includegraphics[width=0.42\textwidth]{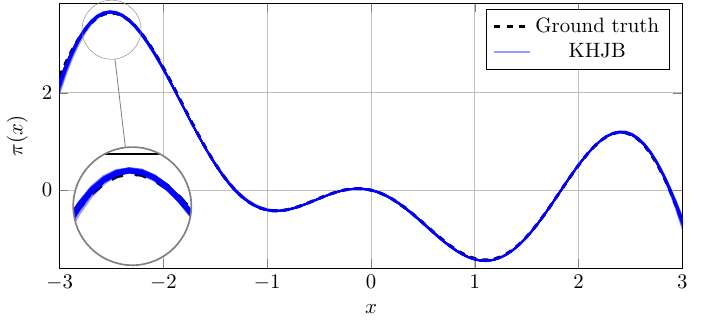}
            \label{fig:subfig3}
        } &
        \subfloat[\texttt{S4}, $\pi^{\star}(x) = x^3 - x \sqrt{1 + x^4}$]{
            \includegraphics[width=0.42\textwidth]{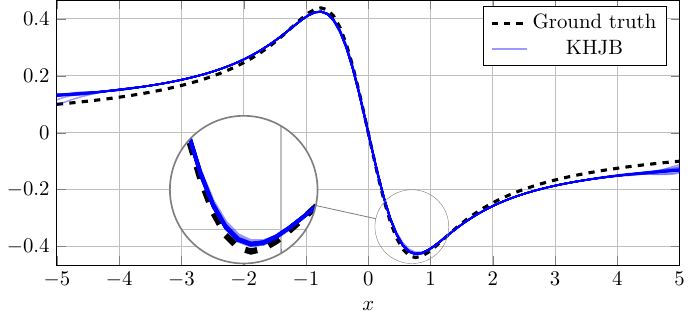}
            \label{fig:subfig4}
        }
    \end{tabular}
    \vspace{-0.3em}
    \caption{The optimal control laws $\pi^\star$ for $\ell(x)=x^2$ and $r(u)=u^2$ are unknown to the algorithm but are used by us as ground truth. Our proposed data-driven Kernel HJB approach recovers these optimal control laws reliably and exhibits little variance across different i.i.d.\ data draws.}
    \label{fig:1dProblems}
\end{figure*}

\subsubsection{1D Benchmark Collection}

Let us consider the open-loop controlled scalar system $\dot{x} =f(x)+G(x) u$ with input $u \in \mathbb{R}$. For $f(x) =-\frac{1}{2} x\left( 1-G(x)^2 \right)$ and stage cost functionals $\stgcost(x)=x^2$, $r(u) =u^2 $, the globally optimal feedback control law is known and
given by ${\pi}^{\star}(x) =-G(x) x$ \citep{Doyle1996}. To obtain control laws with diverse characteristics, we implement three different instances of $G(x)$, listed as Systems~\texttt{S1}, \texttt{S2}, and \texttt{S3} in
Table~\ref{tab:ablation}.  Additionally, we implement another nonlinear system, System~\texttt{S4}, with dynamics $\dot{x}=-x^3+u$, for which the optimal control law has a known closed-form solution  $\pi^{\star}(x)=x^3-x \sqrt{1+x^4}$
\citep{guo2022tutorial}. We generate 50 independent KRR estimations using $n=1000$ training points for Systems~\texttt{S1}, \texttt{S2}, and~\texttt{S3}, and $n=400$ points for System~\texttt{S4}. Algorithm
\ref{alg:ddoc} runs for $H=5000$ steps with sampling time ${h}=10^{-3}$s, and $H=500$ steps with ${h}=10^{-2}$s to ensure value function convergence.  Figure~\ref{fig:1dProblems} compares the control
laws from Algorithm \ref{alg:ddoc} against the known optimal controls.

Table~\ref{tab:rmse1d} reports the mean and standard deviation of the \ac{rmse} $(n^{-1}\sum_{i=1}^n(\widehat{\pi}^\star(\boldsymbol{x}^{(i)}) - \pi^{\star}(\boldsymbol{x}^{(i)}))^2)^{1/2}$ between estimated and optimal control laws. The systems, along with their respective domains and sampling times, are detailed in Table~\ref{tab:ablation}.

\begin{table}[!htbp]\centering
\caption{Definitions of $f$ and $G$ for 1D benchmarks with the respective domain ${\mathbb{X}}_S$ and sampling time ${h}$ of the data.}\label{tab:ablation}

\begin{tabular}{l|cc}
\toprule
System (${\mathbb{X}}_S,{h}$) & $G(x)$ & $f(x)$ \\
\midrule
    \texttt{S1} $\left([-3, 3], 10^{-2}\right)$ & $\sqrt{2}$ & $ \frac{1}{2}x$ \\
    \texttt{S2} $\left([-3, 3], 10^{-3}\right)$ & $\ln{(x^2)}$ & $-\frac{1}{2}x (1 - \ln{(x^2)}^2 )$ \\
    \texttt{S3} $\left([-3, 3], 10^{-3}\right)$ & $\frac{1}{2} + \sin{(2x)}$ & $-\frac{3}{8}x + \frac{1}{2}x \sin{(2x)} +\frac{1}{2}x  \sin{(2x)}^2$\\
    \texttt{S4} $\left([-5, 5], 10^{-2}\right)$ & 1 & $-x^3$
\\
\bottomrule
\end{tabular}
\end{table}

\begin{table}[!htbp]\centering
\caption{1D-benchmarks: \acs{rmse} of the obtained control law versus the ground truth control law, over 100 test data points.}\label{tab:rmse1d}

\begin{tabular}{ll|ll}
\toprule
System & RBF Scale $\sigma$ & \acs{rmse} (mean) & \acs{rmse} (std) \\
\midrule
\texttt{S1} & 1.2 & $1.29 \times 10^{-2}$ & $1.59 \times 10^{-3}$  \\
\texttt{S2} & 1.8 & $1.31 \times 10^{-1}$ & $1.53 \times 10^{-1}$  \\
\texttt{S3} & 2 & $2.17 \times 10^{-2}$ & $8.13 \times 10^{-3}$  \\
\texttt{S4} & 1 & $1.91 \times 10^{-2}$ & $2.94 \times 10^{-4}$
\\
\bottomrule
\end{tabular}
\end{table}

\begin{figure}[!htbp]
    \centering
    \includegraphics[width=0.5\textwidth]{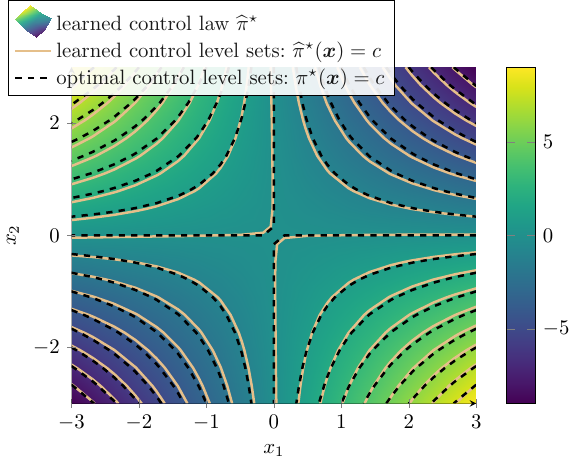}
    \vspace{-0.75em}
    \caption{Contours of the learned optimal controller for the Van der Pol system, compared to the known optimal feedback law $\pi^\star(\boldsymbol{x})=-x_1 x_2$.}
    \label{fig:vdPpolicy}
\end{figure}

\begin{figure}%
    \centering
    \includegraphics[width=0.6\textwidth]{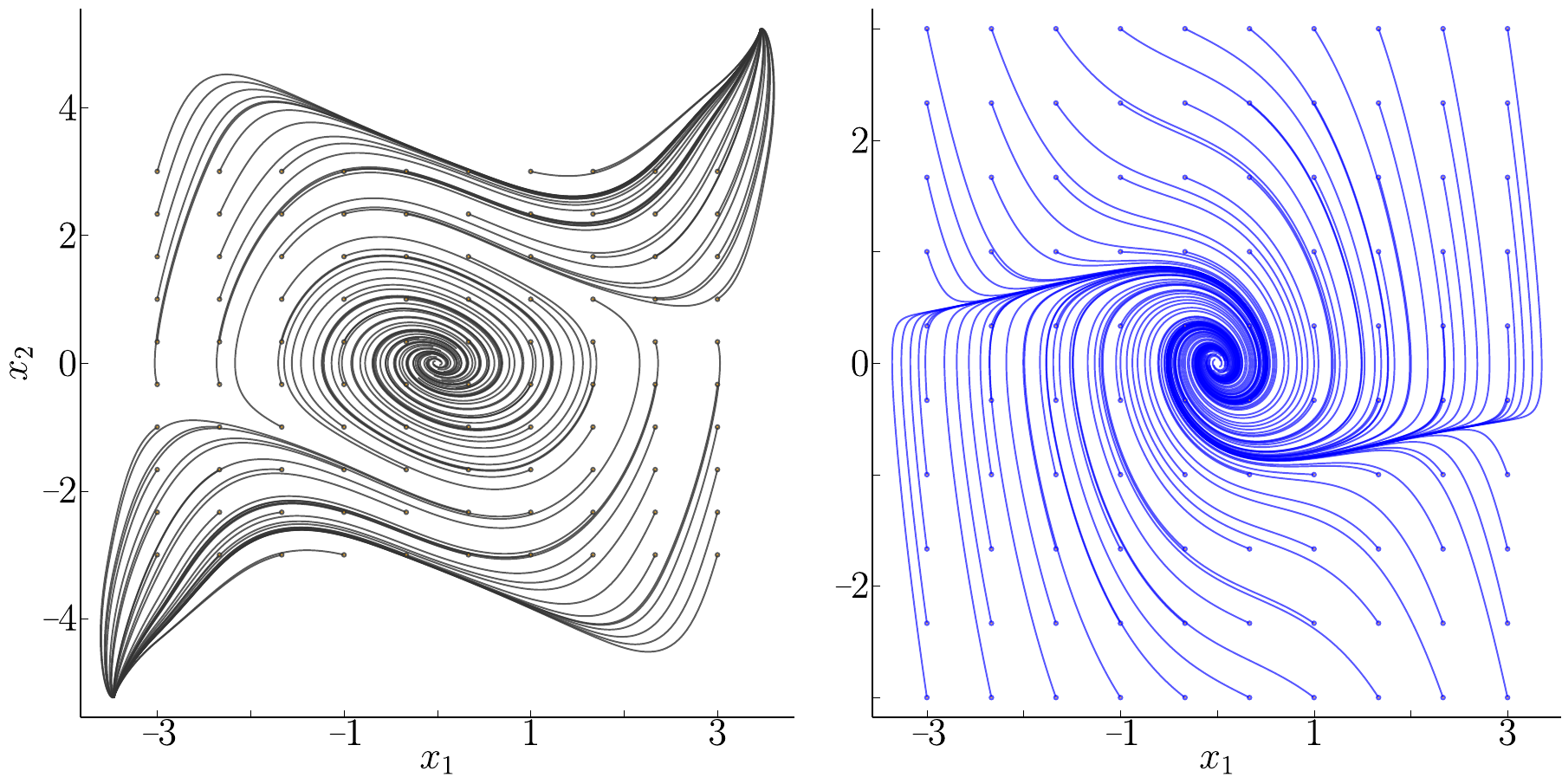}
        \begin{tikzpicture}
    \pgfdeclareplotmark{mystarG}{
    \node[circle,minimum size=1pt, rounded corners=0.1,
          inner sep=0pt,draw=black!50!yellow,solid,fill=black!50!yellow] {};
}
    \pgfdeclareplotmark{mystarR}{
    \node[circle,minimum size=1pt, rounded corners=0.1,
          inner sep=0pt,draw=blue,solid,fill=blue] {};
}
    \begin{axis}[%
    hide axis,
    xmin=0,
    xmax=2,
    ymin=0,
    ymax=0.4,
    domain=0.1:0.1,
    legend style={draw=white!15!black,legend cell align=left,nodes={scale=0.8}},
    legend columns=4,
    width=0.4\columnwidth,
    height=0.4\columnwidth
    ]
    \addlegendimage{gray, thick}
    \addlegendentry{$\texttt{ODEsolve}(\bm{x}_0,0,\texttt{t})$}
    \addlegendimage{blue, thick}
    \addlegendentry{$\texttt{ODEsolve}(\bm{x}_0,\widehat{\pi}(\bm{x}),\texttt{t})$}
    \addlegendimage{mark=mystarG, only marks}
     \addlegendentry{}
    \addlegendimage{mark=mystarR, only marks}
    \addlegendentry{initial state $\bm{x}_0$}
    \end{axis}
\end{tikzpicture}
\vspace{-0.25em}
    \caption{Taming an unstable limit cycling system with linearly uncontrollable origin. \textbf{Left:} Open loop $(u=0)$ of the original system. \textbf{Right:} Closed loop under our learned control law for the stage cost $\ell(\boldsymbol{x})=x^2_2$.}
    \label{fig:VdP_cycle}
\end{figure}

\subsubsection{2D-Unstable oscillator}

Consider a Van der Pol oscillator system
$\dot{x}_1=x_2$, $\dot{x}_2=-x_1-\frac{1}{2} x_2\left(1-{x_1}^2\right)+x_1 u$
with stage cost $\stgcost(\boldsymbol{x})=\frac{1}{2}x_2^2$, $r(u) = \frac{1}{2}u^2$. This system admits an explicit solution to the Hamilton--Jacobi-Bellman equations with globally optimal feedback law $\pi^\star(\boldsymbol{x})=-x_1 x_2$ \citep{Doyle1996}.
The system has a linearly uncontrollable equilibrium at $(0, 0)^\top$ and an unstable limit cycle (Figure~\ref{fig:VdP_cycle}, left). We sample $n=2500$ training points from a grid with $| x_1|  \leq 3$ and $| x_2|  \leq 3$ for KRR estimation using RBF scale $\sigma=20$.
Figure~\ref{fig:vdPpolicy} shows the approximated feedback law closely matches the optimal $\pi^\star(\boldsymbol{x})$, with RMSE $7.24 \times 10^{-2}$ on 900 test points. Under the approximated feedback law, the closed-loop system exhibits an asymptotically stable origin (Figure~\ref{fig:VdP_cycle}, right).

\begin{figure}[ht]
    \centering
    \includegraphics[width=0.5\linewidth]{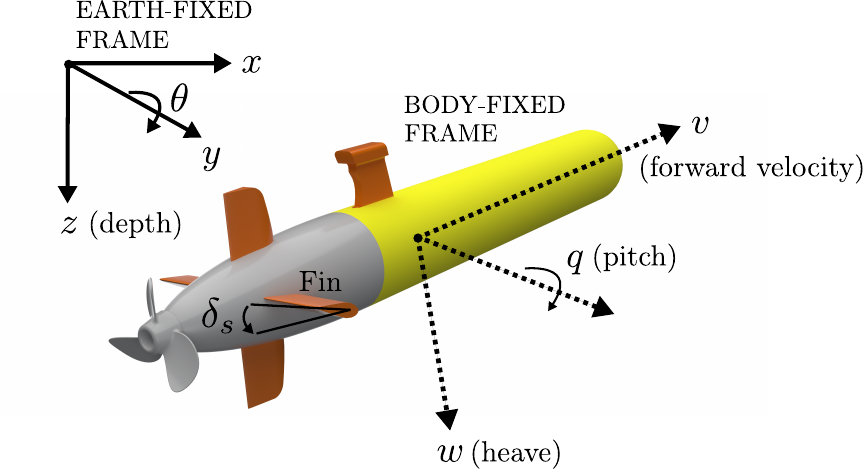}
    \caption{Autonomous Underwater Vehicle (AUV) model.}
    \label{fig:auv_model}
\end{figure}

\subsubsection{4D-Nonlinear {AUV}}

\begin{figure}[!htbp]
    \centering
     \includegraphics[width=1\textwidth]{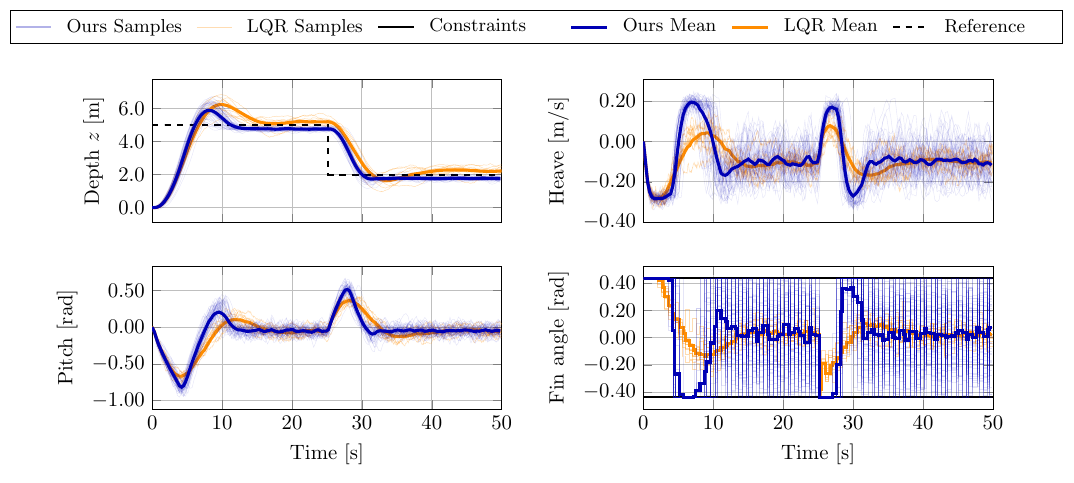}
    \vspace{-1.75em}
    \caption{Tracking performance for two depth references $ z_{\mathrm{ref}} = (5\,\mathrm{m},\,2\,\mathrm{m})$, switching at $ t = (0~\mathrm{s},\, 25~\mathrm{s})$, is evaluated over $M=50$ simulation runs of the closed-loop SDE under the learned control law $ \widehat{\boldsymbol{\pi}}^\star $ with $ \epsilon = 0.001 $ and compared to LQR.}
    \label{fig:AUVrefs}
\end{figure}

We consider depth control of an autonomous underwater vehicle (AUV) depicted in  Figure~\ref{fig:auv_model} with constant forward velocity $v=2~\mathrm{m/s}$ in a vertical plane. The state vector $\boldsymbol{x} = [w, q, z, \theta]^{\tr} \in \mathbb{R}^4$ represents heave velocity
(m/s), pitch rate (rad/s), depth (m), and pitch angle (rad), with control input $\delta_s$ denoting the fin deflection angle. We use the nonlinear model and parameters from \citet{naik2007state}, sampling $n=8000$
i.i.d. training points within $[-0.5, 0.5] \times [-\nicefrac{\pi}{6},\nicefrac{\pi}{6}] \times [0, 4] \times [-\nicefrac{\pi}{6}, \nicefrac{\pi}{6}]$ with random inputs $u \in [-\nicefrac{\pi}{6}, \nicefrac{\pi}{6}]$, collected at a sample rate of $h = 0.5~\mathrm{s}$. Moreover, we identify $\sigma=35$ as the ideal RBF bandwidth and set $\epsilon=0.001$ and $\gamma=10^{-8}$.
For a given depth reference $ z_{\mathrm{ref}} $, the state reference is set as $ \bx_{\mathrm{ref}} = [0, 0, z_{\mathrm{ref}}, 0]^\top $, with quadratic costs $ \ell(\bx) = (\bx - \bx_{\mathrm{ref}})^\top Q (\bx - \bx_{\mathrm{ref}}) $ and $ r(u) = u^\top R u $, where $ Q = \operatorname{diag}(100, 500, 100, 350) + Q_{\mathrm{offdiag}} $ with $ (Q_{\mathrm{offdiag}})_{1,4} = (Q_{\mathrm{offdiag}})_{4,1} = 500 $ and $ R = 50 $. To obtain a
practical control law~\citep{naik2007state}, we constrain the fin-deflection angle. Specifically, we enforce the box-constrained admissible set $\mathcal U=\{u\in\mathbb{R}: | u| \le 0.4363\}$, corresponding to the physical fin-deflection bound $| \delta_s| \le 0.4363~\text{rad}$
(${\approx} \pm25{^\circ}$). In practice, the constrained minimizer is obtained by projecting the unconstrained minimizer onto this interval. The RKHS embedded dynamic programming recursion is executed for $H = 1000$ steps. We then
evaluate the learned control law on the depth-tracking task\footnote{The \emph{same} feedback can track a new reference $z_{\text{ref}}^\prime$ by exploiting the translation invariance of depth in the model of \citet{naik2007state},
redefining $\boldsymbol{\pi}^\star(\bx - [0, 0, z_{\text{ref}}^\prime - z_{\text{ref}}, 0]^\tr)$.}  as shown in  Figure~\ref{fig:AUVrefs} for a simulation time of $T=50~ \mathrm{s}$ and $M=50$ runs.

As the vertical plane model \citep{naik2007state} is made to enable linear control design in state-space for nonlinear systems, we implement a linear quadratic regulator (LQR), which uses the same cost matrices
$Q$ and $R$, and evaluate its performance over the $M$ simulation runs. We compute the accumulated costs $J^{(i)} = \int_0^T \ell\left(\stmbxi(t)\right) + r\left(u^{(i)}(t) \right) \, \mathrm{d}t, \, i \in [M]$, over the trajectories in  Figure~\ref{fig:AUVrefs} for both approaches and find
that our learned control law $\widehat{\boldsymbol{\pi}}$ achieves a mean cost reduction of 5.09\% and a standard deviation reduction of 30.11\% compared to the LQR baseline. The LQR controller, although derived from a linearized model,
benefits from exact knowledge of the system parameters and performs well because it represents an exact linearization at the operating point reached after the transient, due to the translational invariance of depth in the
nonlinear dynamics. In contrast, our control law $\widehat{\boldsymbol{\pi}}$ is learned entirely from data without access to the system model or its parameters, yet achieves lower costs and reduced variability by exploiting the full
nonlinear dynamics.

\section{Conclusion}

This paper has presented a fully data-driven framework for approximating optimal feedback laws in stochastic optimal control problems. This approach \textit{scales linearly with state dimension}, thereby enabling polynomial runtime complexity of both identification as well as control, as summarized in Table~\ref{tab:complexity}. The Kernel HJB recursion and operator model learning offer a promising avenue for addressing high-dimensional and uncertain optimal control problems. This versatility opens a variety of future research directions, including large-scale regression via data-efficient low-rank estimators and feature spaces for further reduction of the computational complexity of the proposed algorithms. Another important direction for future research is to quantify how operator-estimation errors propagate through the dynamic-programming recursion and translate into policy suboptimality. While our framework already naturally accommodates input constraints through the admissible control set, as illustrated by the AUV example, systems requiring stricter state-safety specifications could be addressed by assigning large penalties to unsafe regions in the stage cost, which could be explored in future work.

\newpage
\begingroup
\renewcommand{\bibfont}{\normalsize}
\setlength{\bibhang}{1em}
\setlength{\bibsep}{1em minus 1.5pt}
\bibliographystyle{elsarticle-harv}
\bibliography{references}
\endgroup

\end{document}